\documentclass{article}
\usepackage {geometry}
\usepackage {fancyhdr}
\usepackage {amsmath,amsthm,amssymb}
\usepackage {graphicx}
\usepackage {hyperref}
\usepackage{amsfonts}
\usepackage{tabularx}
\usepackage{enumerate}
\usepackage{undertilde}
\usepackage[normalem]{ulem}

\newtheorem{theorem}{Theorem}

\newtheorem{lemma}[theorem]{Lemma}

\theoremstyle{remark}\newtheorem{definition}[theorem]{Definition}
\theoremstyle{remark}\newtheorem{remark}[theorem]{Remark}

\newcommand{\C}{\mathbb{C}}

\newcommand{\E}{\mathbb{E}}

\newcommand{\Z}{\mathbb{Z}}
\newcommand{\R}{\mathbb{R}}
\renewcommand{\P}{\mathbb{P}}

\newcommand{\GFS}{(C_0^{\infty}(\R^d))'}

\DeclareMathOperator{\GFF}{GFF}
\DeclareMathOperator{\Cov}{Cov}
\DeclareMathOperator{\Var}{Var}
\DeclareMathOperator{\FGF}{FGF}
\DeclareMathOperator{\USF}{USF}

\begin {document}
\author{Xin Sun, Wei Wu}
\title{Uniform Spanning Forests and the bi-Laplacian Gaussian field}

\maketitle

\begin{abstract}
We construct a natural discrete random field on $\mathbb{Z}^{d}$, $d\geq 5$
that converges weakly to the bi-Laplacian Gaussian field in the scaling
limit. The construction is based on assigning i.i.d. Bernoulli random
variables on each component of the uniform spanning forest, thus defines an
associated random function. To our knowledge, this is the first natural
discrete model (besides the discrete bi-Laplacian Gaussian field) that
converges to the bi-Laplacian Gaussian field.

\begin{flushleft}
\textbf{Keywords：}  Uniform spannning forest, bi-Laplacian Gaussian field, moment method
\end{flushleft}
\end{abstract}
\medbreak {\noindent\bf Acknowledgments.} We are very grateful to Scott Sheffield for suggesting this problem, to Gregory Lawler for helping us complete the proof of
Lemma \ref{lemma::correlation} and to Richard Kenyon for discussions.
\section{Introduction} \label{Sec::introduction}
 Uniform spanning forest is an extensively studied combinatorial
object \cite{BLPS}, \cite{Pem91}. The uniform spanning forest measure on $%
\mathbb{Z}^{d}$ can be defined in two equivalent ways: either as the weak
limit of the uniform spanning tree measure on a sequence of finite subgraphs
that exhaust $\mathbb{Z}^{d}$, or as an output of the Wilson's algorithm 
\cite{Wil96}. Detailed descriptions of these constructions are given in Section \ref{Sec:: USF}.

In this paper, we study the following random field associated with the $\USF$
on $\mathbb{Z}^{d}$, $d\geq 5$. It is known that the $\USF$ on $\mathbb{Z}%
^{d}$, $d\geq 5$ has infinitely many tree components a.s. Conditioned on the
configuration of the whole forest $\{T_{i}\}_{i\in \mathbb{N}}$, we assign
i.i.d Bernoulli random variables on each tree $T_{i}$, with probability $1/2$
to be $1$ and $1/2$ to be $-1$. We define a random function (which we call
the spin of the spanning forest) $h_{1}$ from $\mathbb{Z}^{d}$ to $%
\{\pm 1\}$, such that for any $x\in \mathbb{Z}^{d}$, $h_{1}\left( x\right) $
equals the random variable associated with the tree component containing $x$.
This random function is constructed in a similar spirit as the Edward-Sokal
coupling of the FK-Ising model \cite{Gri}.

We would like to study the scaling limit of $h_{1}$. For $ \varepsilon \geq 0 $, consider the lattice $ \varepsilon \Z^d $, let $ h_\varepsilon(x)=\varepsilon^{ \frac{4-d}{2} } h_1(\varepsilon^{-1} x) $, $ \forall x\in \varepsilon \Z^d $.  We extend $ h_\varepsilon $ to $ \R^d $ such that $ h_\varepsilon(y) =  h_\varepsilon(x) $  for $ y\in B_{\varepsilon/2} (x)= (  x-\frac{\varepsilon}{2}, x+\frac{\varepsilon}{2} ]^d $.

It turns out that the limiting field of $ h_\varepsilon $ is a generalized Gaussian field (a random generalized distribution whose integral against any $ C^\infty_0 $ test function is a Gaussian) closely related to bi-Laplacian operator $ \Delta^2 $, which we call the bi-Laplacian Gaussian field. We will give the precise definition of the bi-Laplacian Gaussian field in Section \ref{Sec:: FGF}, here we offer an informal description. Intuitively, a bi-Laplacian Gaussian field is  a generalized Gaussian field $ h $ whose covariance structure is given by $ \Cov[ h(x), h(y)  ]= |x-y|^{4-d} $. The rigorous formulation of this definition of given in Definition \ref{def::def2} of Section \ref{Sec:: FGF}, where we also discuss its relation to bi-Laplacian equations. It is  the analogy of that of Gaussian free field($ \GFF $) to Laplacian equation (for the definition and
properties of Gaussian free field, see the survey \cite{sheffield2007gaussian}). 

Here we point out that the bi-Laplacian Gaussian field fall into a bigger family of Gaussian fields called the fractional Gaussian fields($ \FGF $) which is defined and studied in \cite{FGF}. The relation of  $ \FGF  $ and  fractional Laplacian operator $ (-\Delta)^s $ is analogous to both the Gaussian free field and the bi-Laplacian Gaussian field. Besides $ \GFF  $ and bi-Laplacian free field,  the family of $ \FGF  $ also contains white noise, log-correlated Gaussian field and the fractional Brownian field \cite{adler2007random} (a higher dimensional generalization of fractional Brownian motion).

 The main result of this paper is that $h_\varepsilon $ converge  to  $ h $ as random variables taking values in the space of generalized function. To be precise, we have the following theorem.
\begin{theorem}
\label{thm::main theorem} For any $\varphi \in C_{0}^{\infty }(\mathbb{R}%
^{d})$, $(h_{\varepsilon },\varphi )$ converge to $\sqrt{c_{d}}(h,\varphi )$ in
distribution as $\varepsilon \longrightarrow 0$. The constant $c_{d}$ can be
computed by non-intersecting probability of a simple random walk and two
loop erased random walks, see Lemma \ref{lemma::correlation}.
\end{theorem}

Gaussian fluctuations has been observed and studied for numerous physical
systems. For systems in the critical regime, one expects the spatial or
space-time fluctuation to be Gaussian free field. Typical
examples come from domino tilings \cite{Ken01}, random matrix theory \cite%
{BG13}\cite{rider2006noise} and random growth models \cite{BF08}. In the subcritical regime,
where the correlation decays faster, one expects Gaussian white noise
fluctuations (see the example of edge process of spanning tree models in 
\cite{wu2013USF} ). Our model can be viewed as a natural example in the
supercritical regime.

\cite{sakagawa2003entropic}\cite{kurt2007entropic}\cite{Kur09} study the discrete bi-Laplacian Gaussian field (in physics literature, this is known as the
membrane model) whose continuous counterpart is clearly the bi-Laplacian Gaussian field. Our model can be
viewed as another natural discrete object that converges to the bi-Laplacian
Gaussian field. In one dimensional case, Hammond and Sheffield constructed a
reinforced random walk with long range memory \cite{HS09}, which can be
associated with a spanning forest attached to $\mathbb{Z}$. Our construction
can also be viewed as a higher dimensional analogue of ``forest random walks".

Finally, we remark on universality features of our model. We can replace i.i.d. Bernoulli
random variables by general i.i.d. random variables with mean $0$ and variance $1$, and 
obtain the same scaling limit. The same argument also goes through if we replace $\mathbb{Z}^{d}$
by regular lattices, the constant $c_{d}$ is lattice dependent. See Remark \ref{rmk::universality}.

The strategy of the proof is moment method. Since $(h,\varphi )$ is a
Gaussian random variable, to prove convergence in distribution, we only need
to prove that all the moment of $(h_{\varepsilon },\varphi )$ converge to the
corresponding moments of $(h,\varphi )$. The paper is organized as
follows. In Section \ref{Sec:: preliminary}, we give the necessary
background on uniform spanning forest and the bi-Laplacian Gaussian field. In
Section \ref{Sec:: 2nd moment}, we prove the convergence of second moment.
It involves giving the precise asymptotics of the probability that two
vertices are in the same tree of $\USF$. In Section \ref{Sec:: higher
moments}, we prove the convergence of higher moments. In Section \ref{Sec::
further}, we discuss some further questions.



 \section{Preliminary }\label{Sec:: preliminary}

\subsection{Bi-Laplacian Gaussian field} \label{Sec:: FGF}
In this section, we will give a precise definition of the bi-Laplacian Gaussian field, which is a random variable taking value in the space of generalized functions ( denoted by  $ (C_0^{\infty}(\R^d))' $) .  Or equivalently, a probability distribution on $ \GFS $. For basic facts on generalized function, we refer  to Appendix B, \cite{lax}.

We first review some standard facts on white noise, \cite[text]{kuo1996white}.  White noise is the unique probability distribution on $ \GFS $ such that if $ W $ is a random generalized function with this distribution, then for any $ \varphi \in  C^{\infty}_0(\R^d) $, $ (W,\varphi ) $ is a centred Gaussian variable with variance $ (\varphi, \varphi) $. Here $ (,) $ is the pair of a generalized function and a compact supported smooth function.

Formally speaking, we can say that $ W $ is a Gaussian process whose parameter is in $ \R^d $ and the covariance structure is given by \[  \Cov[(  W(x), W(y )] = \delta(x-y). \] Another natural interpretation is that $ W $ is a standard normal distribution on the Hilbert space $ L^2(\R^d) $.

We give two equivalent definitions of the bi-Laplacian Gaussian field, that only differ by scalar multiplication. We can define the bi-Laplacian Gaussian field for all dimensions in a unified way, as in \cite{FGF}. But  to avoid technical details for $ d\leq 4 $, we only define the field for $ d\geq 5 $, which is sufficient for the purpose of this  paper. From now on, we always assume $ d\geq 5 $.

\begin{definition}\label{def::bi-Laplacian}
Bi-Laplacian Gaussian field is the unique probability distribution on $ \GFS $ such that if $ h $ is a random generalized function with this distribution,  $ \Delta h $ is a white noise on $ \R^d $. Here $ \Delta $ is a well defined operator on $ \GFS $ by integration by part \cite{lax}.
\end{definition}
\begin{definition}\label{def::def2}
Bi-Laplacian Gaussian field is the unique probability distribution on $\GFS  $ such that for any $ \varphi \in C^{\infty }_0 (\R^d)   $, $ (h,\varphi) $ is a centred Gaussian variable and 
\begin{equation}\label{eq::variance}
\Var[ (h,\varphi)  ] = \int \int |x-y|^{4-d} \varphi(x)\varphi(y) dx dy .
\end{equation}
\end{definition}

For this moment, we assume that there is a unique random generalized function satisfying the Definition \ref{def::bi-Laplacian} or \ref{def::def2}, which we will explain later. Now we explain the equivalence of the two definitions. We note that if $ \Delta h= W $. Then $ (h,\Delta^2 f)=(\Delta h, \Delta f)=(W,\Delta f) $ is a centred Gaussian of variance $ (\Delta f, \Delta f)=(f,\Delta^2 f) $. For $ \varphi \in C^{\infty }_0(\R^d) $, we can solve the bi-Laplacian equation 
\begin{equation}\label{eq::bi-laplacian}
 \Delta^2 f = \varphi,
\end{equation}
 for example, using Fourier transform. This is the place our assumption of $ d\geq 5 $ plays a role since otherwise not all functions in $ C^\infty_0(\R^d) $ has bi-Laplacian  inverse. There will be  some extra assumption for $ \varphi $ when $ d\leq 4 $  like in the case of 2 dimensional Gaussian free field in the whole plane\cite{sheffield2007gaussian}. Therefore  the variance of $ (h,\varphi) $ is given by
 \begin{equation}\label{eq::variance}
  (f,\Delta^2 f)=(f,\varphi)
 \end{equation}

The presence of $ \Delta^2 $ is the reason we call the field the bi-Laplacian Gaussian field. In the case of Gaussian free field, we want to solve a Laplacian equation. Here we want to solve a bi-Laplacian equation (\ref{eq::bi-laplacian}). Bi-Laplacian equation is a standard object in potential theory and studied for many years. For information of this equation we refer to \cite{gazzola2010polyharmonic} and references therein. From \cite{gazzola2010polyharmonic}, for $ d\geq 5 $, the fundamental solution of bi-Laplacian equation is $ C_d|x-y|^{4-d} $ and  we can use the fundamental solution  to solve equation (\ref{eq::bi-laplacian}), which is
\begin{equation}\label{eq::convolution}
f(x)= C_d  \int_{\R^d}  |x-y|^{4-d}  \varphi(y)  dy,
\end{equation}
where $ C_d $ is a constant depending on $ d $. From (\ref{eq::variance}) and (\ref{eq::convolution}) we see that Definition \ref{def::bi-Laplacian} and \ref{def::def2} of a bi-Laplacian Gaussian field only differ by a constant $ \sqrt {C_d} $. In Theorem \ref{thm::main theorem} we use Definition \ref{def::def2} as our definition of a bi-Laplacian Gaussian field for convenience.

As mentioned before, the existence and uniqueness in the definition of bi-Lalacian field is not clear as a priori. The rigorous argument for existence and uniqueness is actually the same as in the definition of white noise\cite{kuo1996white}.  Now we sketch the construction of white noise as a random distribution
  following \cite{kuo1996white}. The definition of a bi-Laplacian Gaussian field will follow by a similar argument.
 \begin{definition}[Countably-Hilbert Space]
   Let $V$ be an infinite dimensional vector space over $\C$, and let
   $\{|\cdot|_n\}_{n\geq 1}$ be a collection of inner product norms on
   $V$. Define the metric $d$ on $V$ by 
   \[
   d(u,v) = \sum_{n=1}^\infty 2^{-n}\frac{|u-v|_n}{1+|u-v|_n}, \quad u,v\in
   V,
   \]
   If $V$ is complete with respect to $d$ then $(V,\{|\cdot|_n\}_{n\geq
     1})$ is called a \textit{countably-Hilbert space}.
  \end{definition}
  
  \begin{definition}[Nuclear Spaces]
    Let $V$ be a countably-Hilbert space associated with an increasing
    sequence $\{|\cdot|_n\}_{n\geq 1}$ of norms, that is,
    \[
    |v|_1\leq |v|_2\leq \cdots\leq |v|_n\leq \cdots, \quad \forall v\in V.
    \]
    Let $V_n$ be the completion of $V$ with respect to the norm
    $|\cdot|_n$.  We say that $V$ is a \textit{nuclear space} if for any
    $m$, there exists $n\geq m$ such that the inclusion map $V_m$ into
    $V_n$ is a Hilbert-Schmidt operator, that is, there is an orthonormal
    basis $\{v_k\}$ for $V_m$ such that
    $\sum_{k=1}^\infty|v_k|_n^2<\infty$.
  \end{definition}
  \begin{remark}
    It is well known that $ C^{\infty}_0(\R^d) $ is a nuclear space. For a proof, see
    \cite{kuo1996white}.
  \end{remark}

  If $V$ is a topological vector space, we denote by $V'$ the dual of
    $V$ (that is, the space of continuous linear functionals on $V$). We
    say that a complex-valued function $\varphi$ on $V$ is the
    characteristic function of a probability measure $\nu$ on $V'$ if
    \[
    \varphi(v) = \int_{V'}e^{i(x,v)}\,d\nu(x),\quad \text{for all }v\in
    V.
    \]

  For a proof of the following theorem, see (\cite{kuo1996white}).
  \begin{theorem}[Bochner-Minlos theorem]\label{minlos}

    Let $V$ be a real nuclear space. Then a complex-valued function
    $\varPhi$ on $V$ is the characteristic function of a probability
    measure $\nu$ on $V'$ if and only if $\varPhi(0)=1$, $\varPhi$ is
    continuous, and $\varPhi$ is positive definite, that is,
    
    \[
    \sum_{j,k=1}^nz_j\overline{z_k}\varPhi(v_j-v_k)\geq 0,
    \]
    for all $v_1,\dots,v_n\in V$, and $z_1,\dots,z_n\in \mathbb
    C$. Furthermore, $\varPhi$ determines $\nu$ uniquely. 
  \end{theorem}
  
  White noise will be defined as a Gaussian measure on the space of
  tempered distributions.  To apply Theorem~\ref{minlos}, we first note
  that $C^{\infty}_0 (\R^d)$ is a nuclear space and that the function
  \[
  C(\varphi) = e^{-\frac{1}{2}    (\varphi, \varphi)      },\quad \text{for all
  }\varphi\in C^{\infty}_0 (\R^d),
  \]
  is continuous, positive definite, and satisfies $C(0) = 1$.  Hence
  Theorem~\ref{minlos} implies that there is a unique probability measure
  $\mu$ on $\GFS$ having $C$ as its characteristic
  function. which we define as white noise $W$.  In particular we have the
  relation:
  \[
  \int_{\mathcal S'(\mathbb R^d)}e^{i( x, \varphi)}\,d\mu(x) =
  e^{-\frac{1}{2}(\varphi, \varphi)},\quad \varphi\in C^{\infty}_0 (\R^d),
  \]
  which implies for every $\varphi \in\mathcal C^{\infty}_0(\R^d)$ the random
  variable $(W,\varphi)$ is a mean zero Gaussian with variance
  $(\varphi, \varphi )$.  Given $f,g\in C^{\infty}_0(\R^d)$ we may use polarization to see that 
  \[
    \Cov[(W,f),(W,g)] = (f,g),
  \] We may rewrite the
  above expression as
  \[
    \Cov[(W,f),(W,g)] = \int_{\mathbb R^d}\int_{\mathbb R^d}\delta(x-y) f(x)g(y)\,dx\,dy,
  \]
  and say that $W$ has covariance kernel $\delta(x-y)$.

To show the existence and uniqueness of the bi-Laplacian Gaussian field, we only need to find its characteristic function and apply Theorem \ref{minlos}. From the first definition of Definition \ref{def::bi-Laplacian}, it is easy to see that the characteristic function for a bi-Laplacian Gaussian field is 
\[
C(\varphi) = \frac{1}{\sqrt {2\pi}  } \exp\left(       -\frac{1}{2}     (\Delta^{-1} \varphi, \Delta^{-1} \varphi)      \right).
\]

Here $ (\Delta^{-1} \varphi, \Delta^{-1} \varphi) $ is understood as $ (f,\phi) $ where $ \Delta^2 f=\phi $.

\begin{lemma}
The functional $C(\varphi)$ defined by 
\[
C(\varphi) = \exp\left( -\frac{1}{2} (\Delta^{-1} \varphi, \Delta^{-1} \varphi)   \right),
\]
is a continuous, positive definite functional on $  C^{\infty}_0(\R^d)$ that
satisfies $C(0) = 1$.
\end{lemma}

\begin{proof}
  The continuity (continuity is taken with respect to the norm $ (\Delta^{-1} \varphi, \Delta^{-1} \varphi) ^{  \frac{1}{2}  } $) of $C(\varphi)$  follows from Fourier transform and  the fact that $ |x|^4 $ is integrable in $ \R^d(d\geq 5) $.
  Further the statement $C(0)=1$ is also clear.  All that is left is to check that $C(\varphi)$
  is positive definite.

  Let $\varphi_1,\dots, \varphi_n\in  C^{\infty}_0(\R^d))$ be a set of functions,
  and define $V$ to be the subspace of $C^{\infty}_0(\R^d)$ spanned by
  $\{\varphi_i\}$.  Define $\mu_V$ to be the Gaussian measure on $V$ with
  covariance matrix given by $\Xi_{i,j} = (\Delta^{-1}   \varphi_i, \Delta^{-1}\varphi_j)$ so that its characteristic function is
  \[
  \int_V e^{i(\Delta^{-1} \varphi, \Delta^{-1} y)   }   \,d\mu_V(y) = \frac{1}{\sqrt{2\pi}}  e^{-\frac{1}{2} (\Delta^{-1} \varphi, \Delta^{-1} \varphi)  } = C(\varphi),\quad \varphi\in V.
  \]
  Applying Bochner's theorem for probability measures on $\mathbb R^n$ shows
  us that $C$ is positive definite.
\end{proof}
Now apply Milnos theorem we get the existence and uniqueness of the bi-Laplacian Gaussian field.

\begin{remark}\label{rmk::Def of FGF}
In \cite{FGF}, the authors define the so called fractional Gaussian field in the following way. Formally speaking, the $ d $ dimensional fractional Gaussian field with index $ s $ (denoted by $ \FGF_s^d $)  is given by $ (-\Delta)^{\frac{s}{2}} W $. Thus the bi-Laplacian Gaussian field is $ \FGF_2^d $.
\end{remark}

\subsection{Uniform spanning forest} \label{Sec:: USF}

Here we review some facts about the uniform spanning forest model
(USF) on $\mathbb{Z}^{d}$. Most of the facts extend to general graphs as
well. For more background, we refer the reader to the survey \cite{BLPS}.

Given a finite graph ${\mathcal{G}}\subset \mathbb{Z}^{d}$, the (free)
uniform spanning tree (UST) measure is the probability measure that assign
equal probability to the spanning trees of ${\mathcal{G}}$. When all the
vertices on $\mathbb{Z}^{d}\backslash {\mathcal{G}}$ are contracted to a
single vertex, the corresponding measure is called wired spanning tree.

Uniform spanning forest measure on $\mathbb{Z}^{d}$ is the weak limit of
uniform spanning trees on a sequence of exhausting subsets. Pemantle proved
that the limit of free and wired spanning trees coincide \cite{Pem91}, thus
USF is uniquely defined and has trivial tail.

An alternative way to construct the USF is the Wilson's algorithm \cite{Wil96},
which we now describe. For any path $\mathcal{P}$ in $\mathbb{Z}^{d}$ that
visits each vertex at most finitely many often, the loop erasure of $%
\mathcal{P}$ is constructed as erasing cycles in $\mathcal{P}$ in
chronological order. Fix any ordering $\left( v_{1},v_{2}...\right) $ of
vertices, a growing sequence of forests $\left\{ F_{i}\right\} _{i\in 
\mathbb{N}}$ can be constructed inductively. Let $F_{0}=\emptyset $. Suppose
the forest $F_{i}$ has been generated. Start a simple random walk (SRW) at $%
v_{i+1}$, and stop at the first time it hits $F_{i}$, if it does, and
otherwise let it run indefinitely. $F_{i+1}$ is defined by adding the loop
erasure of this SRW to $F_{i}$ (for $d\geq 3$, SRWs are transient, so the
loop erasure of SRW is well defined a.s.). The algorithm yields $\cup _{i\in 
\mathbb{N}}F_{i}$, it is shown in \cite{BLPS} that its distribution is independent of the ordering of
vertices, and is USF.

Based on Wilson's algorithm and properties of loop erased random walks
(LERW), it is shown in \cite{Pem91} that on $\mathbb{Z}^{d}$, the USF is a
single tree a.s. if $d\leq 4$, and has infinitely many tree components a.s.
when $d\geq 5$. The probablility that two points are in the same tree is the
insection probability of a SRW and a LERW. This will be used in Lemma \ref{lemma::correlation}. Also, when $2\leq d\leq 4$, the
USF has a single topological end a.s. (i.e. removing any vertex disconnect
the tree into two components, one of them is infinite); when $d\geq 5$, each
of the infinitely many trees a.s. has at most two topological ends.



\section{Second moment}\label{Sec:: 2nd moment}
By definition of $ h_\varepsilon $,
\begin{align*}
(h_\varepsilon , \varphi) &= \sum\limits_{x  \in \varepsilon \Z^d  }   \varepsilon^{ \frac{4-d}{2} }  h_1(\varepsilon^{-1} x ) \int_{B_{\varepsilon/2} (x) }  \varphi (y) dy\\
					  &=\sum\limits_{x  \in \varepsilon \Z^d }  \varepsilon^{ \frac{4-d}{2} }  h_1(\varepsilon^{-1} x ) \varphi(x) \varepsilon^{d}+R_\varepsilon(\varphi),
\end{align*}
where the remaining term $\displaystyle\lim_{\varepsilon\rightarrow 0}  R_\varepsilon =0   $ almost surely.

So we only need to show that $ X_\varepsilon=\sum\limits_{x  \in \varepsilon \Z^d }  \varepsilon^{ \frac{4-d}{2} }  h_1(\varepsilon^{-1} x ) \varphi(x) \varepsilon^{d}  $ converge to $ \sqrt{c_{d}}(h,\varphi) $ in distribution. As explained in the introduction, we use the moment method. Since the first moment is just 0, we start from the second moment, which is the focus of this section.

Note that \[  X_\varepsilon =  \sum\limits_{x\in  \Z^d} \varepsilon^{ \frac{4+d}{2} }  h_1(x ) \varphi(\varepsilon x).  \] 
\[ \E[ X^2_\varepsilon ]     = \sum\limits_{x,y\in \Z^d } \varepsilon^{ 4+d } \varphi(\varepsilon x) \varphi(\varepsilon y) \E[h_1(x)h_1(y)].   \]
Let $ p(x,y)=\P[x,y\,\textrm{are in the same tree} ] $, then \[ \E[h_1(x)h_1(y)]=  p(x,y)\times 1+ (1-p(x,y)) \times 0=p(x,y).  \]

As explained in Section \ref{Sec:: USF}, uniform spanning forest can be generated using Wilson algorithm on $ \Z^d $. Therefore from Lemma \ref{lemma::correlation} which we will prove in Section \ref{Sec:: Correlation}, we know that $ \displaystyle\lim_{|x-y|\rightarrow \infty}\frac{p\left( x,y\right) 
}{\left\vert x-y\right\vert ^{4-d}}=c_{d}$. $ c_{d} $ is a constant which we could not explicitly evaluate explicitly because we cannot evaluate the number $ q $ in ( \ref{eq::intersection probability} ) , Section \ref{Sec:: Correlation}.

Since $ \varphi \in C_0^\infty(\R^d) $, by dominate convergence theorem , \[ \displaystyle\lim_{\varepsilon\rightarrow 0} \E[X^2_\varepsilon]  = c_{d}\int\int \varphi(x)\varphi(y)  |x-y|^{4-d} dx dy.   \]
From Section \ref{Sec:: FGF}, we recognize that the RHS of above formula is just  $ c_{d} $  times  the variance of $ (h,\varphi)$ as we defined in formula (\ref{eq::variance}) in Section 
\ref{Sec:: FGF}.

\subsection{Asymptotic correlation}\label{Sec:: Correlation}

In this section we explicitly determine the asymptotics of $p\left( x,y\right)
=p\left( 0,y-x\right) $. This requires to evaluate the non-intersecting
probability of a SRW starts at $y-x$ and a LERW starts at $0$. Using the
bounds for intersection of SRWs, Pemantle showed $p\left( 0,y-x\right)
=O\left( \vert y-x\right\vert ^{4-d} ) $. Here
we show this quantity actually converges in the scaling limit. This requires
a more careful estimate of SRW hitting probabilities. 

\begin{lemma}\label{lemma::correlation} 
Suppose $S^1 $ $S^2 $ are $d(\geq 5) $
dimensional SRWs starting from 0 and $z $. Then $\P ( \hat{S^1} [0,\infty]
\cap S^2[ 0,\infty ] \neq \emptyset ) = c|z|^{4-d} + o(|z|^{4-d})$ as $z \rightarrow
\infty $.
\end{lemma}

\begin{proof}
The proof is suggested by   Lawler \cite{Lawler}.  Let $\rho $ be the first time $S^2 $ hits $\hat{S^1}[0,\infty] $, $\tau $ be
the largest time such that $S^1(\tau ) =S^2( \rho ) $. For $w\in \mathbb{Z}%
^d $, $(j, k)\in \mathbb{N}\times \mathbb{N} $, let $A_{w,j,k} $ be the
event $\{ S^2(\rho)=w, \tau=j, \rho=k \} $. We can see that $\{ \hat{S^1}
[0,\infty] \cap S^2[ 0,\infty ] \neq \emptyset \} =\sum\limits_{w,j,k}
A_{w,j,k} $ almost surely. Let $\bar{S^1} $ be the time reversal of $S^1 $
for $\tau $ to 0, $\bar{S^2} $ be the time reversal of $S^2 $ from $\rho $
to 0, $\bar{S^3} $ be $S^1 $ from $\tau $ to $\infty $. Then 
\begin{equation*}
A_{w, j, k} = \{ \bar{S^1}(j) =0 , \bar{S^2}(k) =z , \hat{\bar{S^1}}
[0,j]\cap \bar{S^3}[1,\infty] =\emptyset, \bar{S^2}[1,k] \cap \{ \hat{ \bar{%
S^1 } }[0,j] \cup \hat{\bar{S^3}} [0,\infty] \} =\emptyset \}. 
\end{equation*}
Thus%
\begin{equation*}
\P (A_{w,j,k}) = \P (\bar{S^1}(j) =0 , \bar{S^2}(k) =z ) \P ( \hat{\bar{S^1}}
[0,j]\cap \bar{S^3}[1,\infty] =\emptyset, \bar{S^2}[1,k] \cap \{ \hat{ \bar{%
S^1 } }[0,j] \cup \hat{\bar{S^3}} [0,\infty] \} =\emptyset| \bar{S^1}(j) =0
, \bar{S^2}(k) =z ). 
\end{equation*}

Now for simplicity of notation we assume that $S^1 , S ^2, S^3 $ are three
independent SRWs starting at $w $. Then 
\begin{equation*}
\P (A_{w,j,k}) = \P (S^1(j) =0 )\P ( S^2(k) =z) \P ( \hat{S^1} [0,j]\cap
S^3[1,\infty] =\emptyset, S^2[1,k] \cap \{ \hat{S^1}[0,j] \cup \hat{S^3}%
[0,\infty] =\emptyset| S^1(j)=0, S^2(k)=z ) . 
\end{equation*}

Let 
\begin{equation}\label{eq::intersection probability}
q=\P ( {\ \hat{S^1} [0,\infty]\cap S^3[1,\infty] =\emptyset, S^2[1,\infty]
\cap \{ \hat{S^1 } [0,\infty] \cup \hat{S^3} [0,\infty] \} } =\emptyset ) 
\end{equation}%
, $G(\cdot,\cdot) $ be the Green function of SRW on $\mathbb{Z}^d $. Now we
show that the non-intersection probability 
\begin{equation*}
\P ( \hat{S^1} [0,\infty] \cap S^2[ 0,\infty ] \neq \emptyset )
=\sum\limits_{w,j,k} \P (A_{w,j,k})\sim q\sum\limits_{w} G(0,w)G(w,z) 
\end{equation*}%
, together with the fact that discrete Green's function converges to the
continuous whole space Green's function \cite{lawler2010random}, therefore $G\left(
z,w\right) =O\left( \left\vert z-w\right\vert ^{2-d}\right) $ for $z,w$
macroscopically apart, this implies Lemma \ref{lemma::correlation}.

To prove the upper bound, we fix small $\varepsilon >0$ and large $R>0$. Let $w$
be in the range of $|w|\geq \varepsilon |z|,|w-z|\geq \varepsilon |z|$ and $j,k$
greater than $|z|^{\frac{3}{2}}$. Let $\sigma _{i}$ be the last time when $%
S^{i}$ hits the ball $B_{R}$ centred at $w$. For fixed $R$ and $w$, on the high probability event that $\sigma _{1}\ll j$ and $\sigma
_{2}\ll k$,  as $|z|\rightarrow \infty $,
the Radon-Nikodym derivative of the joint distribution $\{S^{1}[0,\sigma
_{1}],S^{2}[0,\sigma _{2}],S^{3}[0,\sigma _{3}]\}$ conditioned on that $%
S^{1}(j)=0,S^{2}(k)=z$ w.r.t the original unconditioned one tends to 1. For $%
\left( w,i,j\right) $ satisfy the conditions prescribed, 
\begin{align*}
\P (A_{w,j,k})& \leq \P ^{w}(S^{1}(j)=0)\P ^{w}(S^{2}(k)=z) \\
& \times \P ({\hat{S^{1}}[0,\sigma _{1}]\cap S^{3}[1,\sigma _{3}]=\emptyset
,S^{2}[1,\sigma _{2}]\cap \{\hat{S^{1}}[0,\sigma _{1}]\cup \hat{S^{3}}%
[0,\sigma _{3}]\}}=\emptyset |S^{1}(j)=0,S^{2}(k)=z). \\
& \leq \P ^{0}(S^{1}(j)=w)\P ^{w}(S^{2}(k)=z)(1+\delta _{R,z}) \\
& \times \P ({\hat{S^{1}}[0,\sigma _{1}]\cap S^{3}[1,\sigma _{3}]=\emptyset
,S^{2}[1,\sigma _{2}]\cap \{\hat{S^{1}}[0,\sigma _{1}]\cup \hat{S^{3}}%
[0,\sigma _{3}]\}}=\emptyset ),
\end{align*}%
where $\delta _{R,z}\rightarrow 0$ as $z$ tends to $\infty $ and $R$ fixed.
On the other hand, the typical time for a SRW starting at $w$ to hit $0$ or $%
z$ is $O(\left\vert z\right\vert ^{2})$, thus $\sum\limits_{j<|z|^{\frac{3}{2%
}},k<|z|^{\frac{3}{2}}}\P ^{w}(S^{1}(j)=0)\P ^{w}(S^{2}(k)=z)$ tends to zero
uniformly in $w$ as $z\rightarrow \infty $. Also, when summing over $w\in 
\mathbb{Z}^{d}$, the contribution from $|w|<\varepsilon |z|\,\mathrm{or}%
\,|w-z|<\varepsilon |z|$ is negligible as $\varepsilon \rightarrow 0$. Since%
\begin{equation*}
\lim\limits_{R\rightarrow \infty }\P ({\hat{S^{1}}[0,\sigma _{1}]\cap
S^{3}[1,\sigma _{3}]=\emptyset ,S^{2}[1,\sigma _{2}]\cap \{\hat{S^{1}}%
[0,\sigma _{1}]\cup \hat{S^{3}}[0,\sigma _{3}]\}}=\emptyset )=q.
\end{equation*}%
By summing over $w,i,j$, first taking $z\rightarrow \infty $, then $%
R\rightarrow \infty $ and then $\varepsilon \rightarrow 0$, we know that 
\begin{equation*}
\limsup\limits_{z\rightarrow \infty }\frac{\sum\limits_{w,j,k}\P (A_{w,j,k})%
}{\sum\limits_{w}G(0,w)G(w,z)}\leq q.
\end{equation*}

To show the lower bound, as before we first fix $\varepsilon>0 $, $w$ in the
range $|w|\geq \varepsilon |z|,|w-z|\geq \varepsilon |z|$ and $j,k\geq |z|^{\frac{3%
}{2}}$. When $1\ll R\ll z$ is fixed but large enough, as $z$ tends to $%
\infty $, there is a high probability $p_{R}$ such that the distance between 
$S_{\sigma _{1}}^{1},S_{\sigma _{2}}^{2},S_{\sigma _{3}}^{3}$ is bigger than 
$cR$, where $c$ is a constant independent of $R$ and $p_{R}$ tends to 1 as $R
$ tends to infinity. This is because that as $z\rightarrow \infty $, the $%
S_{\sigma _{1}}^{1},S_{\sigma _{2}}^{2},S_{\sigma _{3}}^{3}$ are close to
three uniform distribution on $\partial B_{R}$ as we argued above. The
probability that $S^{1}[\sigma _{1},\infty ],S^{2}[\sigma _{2},\infty
],S^{3}[\sigma _{3},\infty ]$ have an intersection will tend to zero, as $z$
first goes to $\infty $ and then $R$ goes to $\infty $, which can be seen by
bounding the intersection probabilities explicitly by Green's functions.
Using the asymptotic independence of $S^{1}, S^{2}$ in $B_{R}$ and the event $%
S^{1}(j)=0,S^{2}(k)=z$, we obtain  
\begin{align*}
& \P ( \hat{S^1} [0,j]\cap
S^3[1,\infty] =\emptyset, S^2[1,k] \cap \{ \hat{S^1}[0,j] \cup \hat{S^3}%
[0,\infty] =\emptyset| S^1(j)=0, S^2(k)=z )  \\
& \geq \left( 1-\varepsilon _{R,z}\right) (1-\delta _{R,z})\P ({\hat{S^{1}}%
[0,\sigma _{1}]\cap S^{3}[1,\sigma _{3}]=\emptyset ,S^{2}[1,\sigma _{2}]\cap
\{\hat{S^{1}}[0,\sigma _{1}]\cup \hat{S^{3}}[0,\sigma _{3}]\}}=\emptyset ),
\end{align*}%
where $\delta _{R,z}$ tends to 0 as $z\rightarrow \infty $, and $\varepsilon
_{R,z}\rightarrow 0$ as $z$ first goes to $\infty $ and then $R$ goes to $%
\infty $. Thus 
\begin{equation*}
\liminf\limits_{z\rightarrow \infty }\frac{\sum\limits_{w,j,k}\P (A_{w,j,k})%
}{\sum\limits_{w}G(0,w)G(w,z)}\geq q.
\end{equation*}
\end{proof}



\section{Higher moments} \label{Sec:: higher moments} 

\label{Sec:: higher}

Recall the random field $\{h_{\varepsilon }\}$, defined for any $\varphi
\in C_{0}^{\infty }\left( \mathbb{R}^{d}\right) $ as 
\begin{equation*}
\left( h_{\varepsilon },\varphi \right) =\varepsilon ^{\frac{4+d}{2}%
}\sum_{x\in \varepsilon \mathbb{Z}^{d}}\varphi \left( x\right)
h_{1 }\left( \frac{x}{\varepsilon }\right) +O\left( \varepsilon
\right) \text{.}
\end{equation*}%
Therefore, for $k\geq 3$,%
\begin{eqnarray}
\mathbb{E}\left( \left( h_{\varepsilon },\varphi \right) ^{k}\right) 
&=&\varepsilon ^{\frac{4+d}{2}k}\sum_{x_{1},...,x_{k}\in \varepsilon \mathbb{%
Z}^{d}}\varphi \left( x_{1}\right) ...\varphi \left( x_{k}\right) \mathbb{E}%
\left( h_{1 }\left( \frac{x_{1}}{\varepsilon }\right)
...h_{1 }\left( \frac{x_{k}}{\varepsilon }\right) \right) +O\left(
\varepsilon \right)   \notag \\
&=&\varepsilon ^{\frac{4+d}{2}k}\sum_{\Gamma =\left\{ \gamma _{l}\right\}
}\prod\limits_{l}\sum_{\left\{ x_{m}\right\} _{m\in \gamma _{l}}}\mathbb{E}%
\left( \prod\limits_{m\in \gamma _{l}}\varphi \left( x_{m}\right)
h_{1}\left( \frac{x_{m}}{\varepsilon }\right) \right) +O\left(
\varepsilon \right) .  \label{km}
\end{eqnarray}%
Where in the last equality we group the vertices in terms of components of
the uniform spanning forest: we sum over all the partitions $\Gamma $ of
the index set $\left\{ 1,...,k\right\} $, $h_{1}$ at vertices
belong to different components of the forest are independent.

We claim the following Wick's formula holds in the limit:%
\begin{equation*}
\lim_{\varepsilon \rightarrow 0}\mathbb{E}\left( \left( h_{\varepsilon
},\varphi \right) ^{k}\right) =\left\{ 
\begin{array}{cc}
\left( k-1\right) !!\left( \lim_{\varepsilon \rightarrow 0}\mathbb{E}\left(
\left( h_{\varepsilon },\varphi \right) ^{2}\right) \right) ^{k/2} & k\text{
even} \\ 
0 & k\text{ odd}%
\end{array}%
\right. ,
\end{equation*}%
It therefore uniquely identify the distribution of $\lim_{\varepsilon
\rightarrow 0}\left( h_{\varepsilon },\varphi \right) $ to be Gaussian. By
the covariance structure given is Section \ref{Sec:: 2nd moment}, we
complete the proof that $h_{1 }$ converges weakly to $h$.

When $k$ is odd, at least one of the $\gamma _{l}$ contains odd number of
elements, and therefore $\mathbb{E}\left( \prod\limits_{m\in \gamma
_{l}}h_{1}\left(  \frac{x_{m}}{\varepsilon}     \right) \right) =0$. The independence of $%
h_{1}$ at different components implies $\mathbb{E}\left( \left(
h_{\varepsilon },\varphi \right) ^{k}\right) =0$.

When $k$ is even, the non-vanishing contribution only comes from partitions
such that each $\gamma _{l}$ contains even number of elements. By (\ref{km}%
), it suffices to show that the contribution from those $\left\{ \gamma
_{l}\right\} $, with some $\left\vert \gamma _{l}\right\vert \geq 4$ is
negligible in the limit. We claim: 
\begin{equation*}
\mathbb{E}\left( \prod\limits_{m=1}^{2l}h_{1}\left( \frac{x_{m}}{\varepsilon} \right)
\right) =O\left( \varepsilon ^{\left( d-4\right) \left( 2l-1\right) }\right)
.
\end{equation*}%
And therefore, the contribution from the partition with a cycle of length $2l
$ is 
\begin{equation*}
\varepsilon ^{\frac{4+d}{2}2l}\sum_{x_{1},...,x_{2l}}\mathbb{E}\left(
\prod\limits_{m=1}^{2l}\varphi \left( x_{m}\right) h_{1}\left( 
\frac{x_{m}}{\varepsilon }\right) \right) \leq O\left( \varepsilon ^{\frac{%
4+d}{2}2l}\varepsilon ^{-2dl}\varepsilon ^{\left( d-4\right) \left(
2l-1\right) }\right) =O\left( \varepsilon ^{\left( d-4\right) \left(
l-1\right) }\right) ,
\end{equation*}%
which vanishes for $l\geq 2$.

Note that $\mathbb{E(}\prod\limits_{m=1}^{2l}h_{1 }\left(
\frac{x_{m}}{\varepsilon} \right) )$ is the probability that $\frac{x_{1}}{\varepsilon},...,\frac{x_{2l}}{\varepsilon}$ belong to the
same tree component. This can be computed in terms of intersection
probability of LERWs by Wilson's algorithm (see Section \ref{Sec:: USF}). It
is given by the probability of the following event: start a LERW from $x_{1}/\varepsilon$%
, and run indefinitely; then for $m=2,...,2l$, start a SRW from $x_{m}/\varepsilon$ that
eventually hit the union of the $m-1$ walks, then stopped, and add its loop
erasure to the union of the $m-1$ walks. Since LERW is a subset of SRWs, the
above quantity is bounded by the corresponding intersecting events of SRWs. 
The probability of each of such events can be bounded in a simple way. We prove it in detail
for one example, the others are similar.
For instance, let $A(x_{1},...,x_{2l})$ denote the event, that the SRW starting at $x_{2}/\varepsilon$
hits the SRW starting at $x_{1}/\varepsilon$, the SRW starting at $x_{3}/\varepsilon$
hits the SRW starting at $x_{2}/\varepsilon$, and so on. Then
\begin{eqnarray*}
&&\mathbb{P}\left( A(x_{1},...,x_{2l}) \right)  \\
&\leq &\sum_{w_{1},...,w_{2l-1}\in \mathbb{Z}^{d}}\mathbb{P}%
\left( SRW_{x_{1}/\varepsilon }\text{ hit }w_{1}\text{; }SRW_{x_{2}/%
\varepsilon }\text{ hit }w_{1},w_{2};...;SRW_{x_{2l-1}/\varepsilon }\text{
hit }w_{2l-2},w_{2l-1};SRW_{x_{2l}/\varepsilon }\text{ hit }w_{2l-1}\right) 
\\
&\leq &\sum_{w_{1},...,w_{2l-1}\in \mathbb{Z}^{d}}G\left(
x_{1}/\varepsilon ,w_{1}\right) G\left( x_{2},w_{1}\right) G\left(
x_{2}/\varepsilon ,w_{2}\right) ...G\left( x_{2l}/\varepsilon
,w_{2l-1}\right)  \\
&=&\left( \sum_{w_{1}\in \varepsilon \mathbb{Z}^{d}}G\left( x_{1}/\varepsilon
,w_{1}/\varepsilon \right) G\left( x_{2}/\varepsilon ,w_{1}/\varepsilon
\right) \right) ...\left( \sum_{w_{2l-1}\in \varepsilon \mathbb{Z}%
^{d}}G\left( x_{2l-1}/\varepsilon ,w_{2l-1}/\varepsilon \right) G\left(
x_{2l}/\varepsilon ,w_{2l-1}/\varepsilon \right) \right)  \\
&=&O\left( \varepsilon ^{\left( d-4\right) \left( 2l-1\right) }\right) ,
\end{eqnarray*}%
where the second inequality follows from the fact that the probability of a
SRW hitting a point is bounded by the expected hitting time, which is given
by the lattice Green's function. The last inequality follows from the
Green's function asymptotics $G\left( x/\varepsilon ,w/\varepsilon \right)
=O\left( \varepsilon ^{d-2}\right) $ \cite{lawler2010random}. 
Since $\mathbb{E}\left( \prod\limits_{m=1}^{2l}h_{1 }\left(
x_{m}/\varepsilon \right) \right)$ is a sum of finitely many such probabilities,
it is at most $O\left( \varepsilon ^{\left( d-4\right) \left( 2l-1\right) }\right)$. And the proof is complete.

\begin{remark}\label{rmk::universality}
From the argument in Section \ref{Sec:: 2nd moment} and \ref{Sec:: higher moments}, we can see that the proof does not require many special properties of 
Bernoulli random variables. What we need is that the sequence of i.i.d random variables have mean 0, variance 1, and all finite moments. Moreover, on other regular lattices, since the Green's function has the same asymptotic decay rate (because the SRW still converges to Brownian motions), our result also holds for  uniform spanning forest on other regular lattices. In this sense, Theorem \ref{thm::main theorem} is
universal.
\end{remark}


\section{Further questions}\label{Sec:: further}
\begin{enumerate}

\item Bi-Laplacian Gaussian field is conformally invariant in four dimension.
Are there any discrete random fields on $\mathbb{Z}^{4}$ that scale to
some bi-Laplacian Gaussian field?

\item What geometric properties of uniform spanning forest can be inferred from
the bi-Laplacian Gaussian field?

\item If one introduces short range interactions between the spins on different trees, do one obtain the same scaling limit?

\end{enumerate}


\bibliographystyle{plain}
\bibliography{USFFGFRef}

\bigskip

\filbreak
\begingroup
\small
\parindent=0pt

\bigskip
\vtop{
\hsize=5.3in
Xin Sun\\
Department of Mathematics\\
Massachusetts Institute of Technology\\
Cambridge, MA, USA \\
xinsun89@math.mit.edu

} 

\bigskip
\vtop{
\hsize=5.3in
Wei Wu\\
Division of Applied Mathematics\\
Brown University\\
Providence, Rhode Island, USA \\
wei$\_ $wu@brown.edu
}
\endgroup \filbreak

\end{document}